%% file: MCGLochNess.tex
\documentclass{amsart}

\input{header}

\input{defn}

\title{A note on subgroups of the Loch Ness Monster Surface's mapping class group}

\author{Yannick Krifka}
\email{krifka@mpim-bonn.mpg.de}
\address{Max-Planck-Institut für Mathematik,
Vivatsgasse 7,
53111 Bonn,
Deutschland}

\author{Davide Spriano}
\email{spriano@maths.ox.ac.uk}
\address{Mathematical Institute,
University of Oxford,
Andrew Wiles Building,
Radcliffe Observatory Quarter,
Woodstock Road,
Oxford,
OX2 6GG}

\date{\today}

\begin{document}
    \input{abstract}

    \maketitle
    
    \input{introduction}
    
    \input{proof}

    \input{acknowledgements}
    
     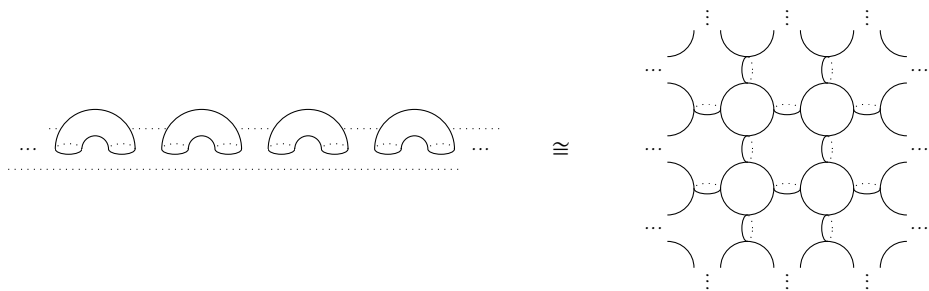
\begin{figure}[H]
     \centering
     \begin{tikzpicture}[scale=0.7]
         \draw[dotted] (0,0,1) -- (8.5,0,1);
         \draw[dotted] (0,0,-1) -- (8.5,0,-1);
         \node[transform shape] at (0,0,0) {$\ldots$};
         \node[transform shape] at (8.5,0,0) {$\ldots$};
         \foreach \x in {0.5, 2.5, 4.5, 6.5}
             \pic[scale=0.5, transform shape] at (\x,0) {handle};
            
          \node at (10,0) {$\cong$};
        
         \begin{scope}[xshift = 12cm, yshift = - 1.75cm, scale = 0.5]
             \foreach \x in {0,3,6}
                 \foreach \y in {0,3,6}
                     \pic[transform shape] at (\x, \y) {xpiece};
                    
             \foreach \y in {0.5,3.5,6.5} {
                 \node[scale = 0.7] at (-0.5,\y) {$\ldots$};
                 \node[scale = 0.7] at (9.5,\y) {$\ldots$};
                 }
            
             \foreach \x in {1.5,4.5,7.5} {
                 \node[scale = 0.7, rotate=90] at (\x,-1.5) {$\ldots$};
                 \node[scale = 0.7, rotate=90] at (\x,8.5) {$\ldots$};
             }
            
             \foreach \x in {3,6}
                 \foreach \y in {0,3,6}
                     \pic[yshift = 0.5cm,transform shape] at (\x, \y) {surfcircv};
            
             \foreach \x in {1.5,4.5,7.5}
                 \foreach \y in {2,5}
                     \pic[transform shape] at (\x, \y) {surfcirch};
         \end{scope}
     \end{tikzpicture}
     \caption{The Loch Ness monster surface $L$ of \cite{sullivan, ghys}, also known as the infinite prison window.}
     \label{fig:lochnessmonster}
 \end{figure}
    
    \printbibliography

\end{document}

%% file: header.tex
\usepackage[utf8]{inputenc}
\usepackage{amsmath}
\usepackage{amsthm}
\usepackage{amssymb}
\usepackage{dsfont}
\usepackage{mathtools}
\usepackage{tikz-cd}
\usepackage{tikz}
\usetikzlibrary{decorations.pathreplacing}
\usetikzlibrary{decorations.pathmorphing}
\usetikzlibrary{decorations.markings}
\usetikzlibrary{arrows.meta}
\usepackage{enumitem}
\usepackage{import}
\usepackage{subcaption}
\usepackage{float}
\usepackage[section]{placeins}

\usepackage{kpfonts}

\usepackage{tabu}
\usepackage{geometry} \usepackage{xcolor}

	\definecolor{todocolor}{HTML}{73c7ff}
	\definecolor{lightblue}{HTML}{9dc4f5}
	\definecolor{lightred}{HTML}{f59d9d}
	
	\definecolor{darkblue}{HTML}{0502b3}
	\definecolor{darkgreen}{HTML}{0eb302}
	\definecolor{darkorange}{HTML}{e08e00}
	\definecolor{darkred}{HTML}{b30202}
	
	\definecolor{linkcol}{HTML}{004f8c}
	\definecolor{citecol}{HTML}{173775}

\usepackage[color=todocolor!40]{todonotes}
\usepackage[onehalfspacing]{setspace}

\usepackage{blindtext}

\usepackage[
giveninits=true,
backend=bibtex,
style=alphabetic,
sorting=nyt,
maxbibnames = 50,
hyperref=false,
url=false,
date=year
]{biblatex}

\usepackage[hyperfootnotes=false]{hyperref}

\hypersetup{
	breaklinks=true,
	colorlinks=true,
	linktocpage=true,
	linkcolor=linkcol,
	anchorcolor = linkcol,
	citecolor=blue,
	urlcolor=blue
}

\AtEveryBibitem{\clearlist{language}}
\allowdisplaybreaks

%% file: defn.tex
\newcommand{\ZZ}{\mathbb{Z}}

\renewcommand{\epsilon}{\varepsilon}

\renewcommand{\i}[1]{#1^{-1}}

\newcommand{\Homeo}{\operatorname{Homeo}}

\renewcommand{\phi}{\varphi}
\newcommand{\MCG}{\operatorname{MCG}}

\newcommand{\Cay}{\operatorname{Cay}}

\renewcommand{\d}{\,\,d}

\renewcommand{\O}{\mathcal{O}}

\tikzset{
	hole/.pic = {
		\draw (-0.5,0) to[bend right] (0.5,0);
		\draw (-0.3,-0.1) to[bend left] (0.3,-0.1);
	}
}

\tikzset{
	puncture/.pic = {
		\fill[white] (0,0) circle (0.1);
		\fill[black] (0,0) circle (0.05);
	}
}

\tikzset{
	marking/.pic = {
		\draw[thick] (-0.05,-0.05) to (0.05,0.05);
		\draw[thick] (-0.05,0.05) to (0.05,-0.05);
	}
}

\tikzset{
	emptypuncture/.pic = {
		\fill[white] (0,0) circle (0.1);
	}
}

\tikzset{
    handle/.pic = {
        \draw[fill, white] (0,0,0) to[out=90,in=180] (1.5,1.5,0) to [out=0,in=90] (3,0,0)
            to (2,0,0) to[out=90, in=0] (1.5,0.5,0) to[out=180,in=90] (1,0,0) to (0,0,0);
        \draw (0,0,0) to[out=90,in=180] (1.5,1.5,0) to [out=0,in=90] (3,0,0);
        \draw (1,0,0) to[out=90,in=180] (1.5,0.5,0) to [out=0,in=90] (2,0,0);
        \begin{scope}
            \clip (0,0) rectangle (1,0.2);
            \draw[dotted] (0.5,0) ellipse (0.5 and 0.2);
        \end{scope}
        \begin{scope}
            \clip (0,-0.2) rectangle (1,0);
            \draw (0.5,0) ellipse (0.5 and 0.2);
        \end{scope}
        
        \begin{scope}
            \clip (2,0) rectangle (3,0.2);
            \draw[dotted] (2.5,0) ellipse (0.5 and 0.2);
        \end{scope}
        \begin{scope}
            \clip (2,-0.2) rectangle (3,0);
            \draw (2.5,0) ellipse (0.5 and 0.2);
        \end{scope}
        
    }
}

\tikzset{
    xpiece/.pic = {
        \draw (0,1) to[out=0,in=-90] (1,2);
        \draw (2,2) to[out=-90,in=180] (3,1);
        \draw (3,0) to[out=180,in=90] (2,-1);
        \draw (1,-1) to[out=90,in=0] (0,0);
    }
}

\tikzset{
    surfcircv/.pic = {
        \begin{scope}
            \clip (0,-0.5) rectangle (0.2,0.5);
            \draw[dotted] (0,0) ellipse (0.2 and 0.5);
        \end{scope}
        \begin{scope}
            \clip (-0.5,-0.5) rectangle (0,0.5);
            \draw (0,0) ellipse (0.2 and 0.5);
        \end{scope}
    }
}

\tikzset{
    surfcirch/.pic = {
        \begin{scope}
            \clip (-0.5,0) rectangle (0.5,0.2);
            \draw[dotted] (0,0) ellipse (0.5 and 0.2);
        \end{scope}
        \begin{scope}
            \clip (-0.5,-0.5) rectangle (0.5,0);
            \draw (0,0) ellipse (0.5 and 0.2);
        \end{scope}
    }
}

\newcommand{\p}[1]{\left( #1 \right)}

\theoremstyle{plain}

\newtheorem*{thm*}{Theorem}
\newtheorem{thm}{Theorem}[section]

\theoremstyle{definition}

\theoremstyle{remark}
\newtheorem{remark}[thm]{Remark}

\definecolor{inkscapepurple}{rgb}{0.50196078,0,0.50196078}
\definecolor{inkscapeblue}{rgb}{0,0,1}

\expandafter\def\expandafter\UrlBreaks\expandafter{\UrlBreaks
	\do\a\do\b\do\c\do\d\do\e\do\f\do\g\do\h\do\i\do\j%
	\do\k\do\l\do\m\do\n\do\o\do\p\do\q\do\r\do\s\do\t%
	\do\u\do\v\do\w\do\x\do\y\do\z\do\A\do\B\do\C\do\D%
	\do\E\do\F\do\G\do\H\do\I\do\J\do\K\do\L\do\M\do\N%
	\do\O\do\P\do\Q\do\R\do\S\do\T\do\U\do\V\do\W\do\X%
	\do\Y\do\Z}

%% file: abstract.tex
\begin{abstract}
    In this short note we give an elementary proof of the fact that every countable group is a subgroup of the mapping class group of the Loch Ness monster surface.
\end{abstract}

%% file: introduction.tex
\section{Introduction}





The following theorem follows from a deep result of Aougab--Patel--Vlamis \cite{aougab_patel_vlamis} about hyperbolic isometry groups of infinite type surfaces:

\begin{thm} \label{thm:main}
    Every countable group $G$ is a subgroup of $\MCG(L)$ where $L$ denotes the Loch Ness monster surface.
\end{thm}

The purpose of this short note is to provide an elementary proof of this fact.

%% file: proof.tex
\section{Proof of Theorem \ref{thm:main}}

    By the Higman--Neumann--Neumann embedding theorem \cite[Theorem IV]{higman} any countable group $G$ is a subgroup of a finitely generated group $H'$. Thus, $G$ is also a subgroup of the product $H \coloneqq H' \times \ZZ \times \ZZ$. 
    Let $S'$ denote a finite generating set for $H'$. Then $S = \{ (s',0,0) \, | \, s' \in S'\} \cup \{ (e,1,0), (e,0,1)\} \subseteq H' \times \ZZ \times \ZZ$ is a generating set of $H$. Observe that the Cayley graph $\Gamma \coloneqq \Cay_S(H)$ is one-ended  and has degree at least $4$. Therefore, the surface $\Sigma$ that is obtained by thickening the Cayley graph $\Gamma$ of $H$ is one-ended, too. Moreover, one may check that $\Sigma$ has infinite genus. By the classification of topological surfaces of infinite type \cite{kerekjarto, richards} $\Sigma$ is thus homeomorphic to the Loch Ness monster surface $L$; see Figure \ref{fig:lochnessmonster}.
    
    As a subgroup of $H$, the group $G$ acts freely on $\Gamma$ via graph automorphisms. Thus, we can choose an invariant hyperbolic metric on $\Sigma$ so that this action gives rise to a free action of $G$ on $\Sigma \cong L$ by orientation preserving isometries. As every non-trivial isometry is a homotopically non-trivial homeomorphism \cite{Norris}, we obtain an injection $ G \hookrightarrow \Homeo_+(L)/\Homeo_+^\circ(L) = \MCG(L)$.\qed

    

    \begin{remark}
    The same result holds for the blooming Cantor tree surface instead of the Loch Ness monster surface. This is obtained by replacing the group $H$ with the group $K = H' \ast \ZZ  \ast \ZZ/2\ZZ$ in the proof. 
    \end{remark}


%% file: acknowledgements.tex
\textsc{Acknowledgements:} We are grateful to Nicholas Vlamis for useful feedback on an early draft of the paper, and for pointing out that our method carries over to the blooming Cantor tree surface. 